\documentclass[11pt]{amsart}
    \usepackage{parskip}
    \usepackage[margin=1in]{geometry}
    \usepackage{amsmath,amssymb,amsthm}
    \usepackage[expansion=false]{microtype}
    \usepackage[colorlinks=true,linkcolor=blue,citecolor=blue,urlcolor=blue]{hyperref}
    \usepackage{xcolor}                                   
    \newtheorem{theorem}{Theorem}
    \newtheorem{proposition}[theorem]{Proposition}
    \newtheorem{corollary}[theorem]{Corollary}            
    \newtheorem{lemma}[theorem]{Lemma}                    
    \theoremstyle{remark}

    \newcommand{\Q}{\mathbf{Q}}
    \newcommand{\Z}{\mathbf{Z}}
    
    \newcommand{\cC}{\mathcal C}
\newcommand{\cK}{\mathcal K}
\newcommand{\cD}{\mathcal D}
\newcommand{\cR}{\mathcal R}
\DeclareMathOperator{\Vol}{Vol}

    \title{Consequences of ``No Siegel zeros''}
    \author{Andrew Granville}
    \address{ }
    \email{}
    \subjclass[2020]{Primary 11N13; Secondary 11N36}
    \date{}
    \begin{document}
    \begin{abstract} 
    We give a more precise and arguably easier version of OpenAI's proof that there are no Siegel zeros. We give uniform bounds on  
    $ \sum_{ p\leq x } \frac{ (d/p) \log p}p$ and therefore that $ | \frac{L'}{L}(1,(d/\cdot))  | \leq \ \log |d|+\log\log |d| +O(1)$. We deduce that 
    if $L(\beta,(d/\cdot))=0$ with $\beta\in \mathbb R$    then $ \beta\leq 1-\frac{c}{\log |d|}$ for any  $c<1$; and show, in an elementary way, that
    $L(1,(d/\cdot))\ \ge\ \frac{c}{\log|d|} $  for any  $c<1/5$. Finally we develop a preprint of Pollack \cite{poll} using ideas of Matom\"aki,   Merikoski, and   Ter\"av\"ainen \cite{mmt} to give an elementary explicit Linnik's theorem; that 
      there is a prime $p\equiv a\pmod Q$ with $p\le Q^{138}$.  Claude helped me with some of the more technical aspects.
    \end{abstract}
    \maketitle
    \vspace{-0.2in}
    \section{Introduction}  This is a moment of great change for analytic number theory. Our main goal in this article is to explore the ``obvious'' mainstream consequences of OpenAI's recent proof that there are no Siegel zeros, in preparation for teaching analytic number theory next semester   (and sharing these notes for others in the same predicament).

    \subsection{$L(s,(d/\cdot))$ close to $1$, and $(L'/L)(1,(d/\cdot))$.}
      Throughout, $d\neq 1$ is a fundamental discriminant, $\chi_d=(d/\cdot)$, and $p,q$ denote primes.

 Our first goal will be to give bounds for 
$\sum_{ p\leq x } \frac{(d/p) \log p}p $ which leads to suitable bounds on real zeros of $L(s,(d/\cdot))$ (so that there are no ``Siegel zeros'').
The initial bounds in our proof are completely elementary, but by combining these with the explicit formula, the classical zero-free region, Dirichlet's class number formula and a simple bootstrapping method, we prove the following:

\begin{corollary}\label{thm:L'/L}
Let $d$ be a fundamental discriminant.  Then
 \[
 \bigg| \frac{L'}{L}(1,(d/\cdot)) \bigg| \leq \ \log |d|+ \log\log |d|+O(1).
 \]
 Here the ``$O(1)$'' is   effective.
 \end{corollary}

 GRH gives the stronger result
  \[
 \bigg| \frac{L'}{L}(1,(d/\cdot)) \bigg| \ll \log\log |d|.
 \]

\begin{corollary}\label{cor:L2}
If  $L(\beta,(d/\cdot))=0$ with $\beta\in \mathbb R$  for a  fundamental discriminant $d$ then
\[
 \beta\leq 1-\frac 1{ \log |d|+3\log\log |d| + 2\log\log\log |d|+O(1)}.
\]
\end{corollary}

GRH gives the stronger result $\beta\leq \frac 12$.

\begin{corollary}\label{cor:L1}
For every fundamental discriminant $d$ we have
\[
L(1,(d/\cdot))\ \ge\ \frac{1}{5\log|d|}  \bigg( 1-O \bigg( \frac{\log\log |d|}{\log |d|}\bigg)\bigg).
\]
\end{corollary}

GRH gives the stronger result $L(1,(d/\cdot))\ \ge\  \frac{e^{-\gamma}\pi^2}{12\log\log|d|+O(1)}$.

\subsection{Linnik's theorem}
Linnik's theorem is a well-known difficult-to-prove result in analytic number theory stating that there exists a constant $L$ such that if
$(a,Q)=1$ then there is a prime $p\equiv a\pmod Q$ with $p\ll Q^L$.  Until recently the best unconditional result was $L=5$, due to Xylouris \cite{Xyl} (improving his earlier $L=5.18$ \cite{Xyl11}), while under GRH one can take any $L>2$ with little difficulty.

\subsubsection{No Siegel zeros  and Linnik} 
In this article we are concerned with elementary proofs, particularly those using the recent breakthrough on Siegel zeros.
Only recently Matom\"aki,   Merikoski, and   Ter\"av\"ainen \cite{mmt} gave a beautiful elementary proof based on using the small sieve (inspired by work of Friedlander and Iwaniec \cite{fi23}). 
Pollack \cite{poll}  used Selberg type identities to show if there are no small primes in an arithmetic progression then there are few integers that are the product of three small primes in the same arithmetic progression and then used  this breakthrough on Siegel zeros   to prove Linnik's theorem. By a seemingly very different method Matom\"aki,   Merikoski, and   Ter\"av\"ainen \cite{mmt} also study products of three primes instead of the primes themselves, so we have adapted 
the proof in \cite{mmt} to the information that Pollack uses and we obtain the following explicit result:

\begin{theorem}\label{thm:primeslinnik}
For $Q$ sufficiently large, for $x= Q^{138}$ and every $(a,Q)=1$,
\[
 \sum_{\substack{p\le x\\ p\equiv a\ (\mathrm{mod}\ Q)}}\frac{\log p}{p}\ \ge\ \frac{\log x}{3600\,\phi(Q)}.
\]
In particular there is a prime $p\equiv a\pmod Q$ with $p\le Q^{138}$.
\end{theorem}

  Matom\"aki,   Merikoski, and   Ter\"av\"ainen \cite{mmt} obtained the exponent $350$; so our exponent is only a little smaller, but the proof is somewhat simpler.

\subsection*{Thanks} Thanks are due to Dimitris Koukoulopoulos and Sun-Kai Leung for helpful suggestions.

\subsection*{Use of AI} The author used Claude to help develop the technical details, but all of  the ideas come from the references, the standard literature, Dimitris Koukoulopoulos and Sun-Kai Leung , as well as the author.  The author reviewed and verified the AI-assisted material and takes full responsibility for the mathematical content of the manuscript.

\section{Bounds on partial sums of $ (L'/L)(1,\chi)$ from first principles}

We recall the well-known estimate
\begin{equation} \label{eq: Trivial}
\sum_{p\leq x} \frac{ \log p}p =\log x+O(1).
\end{equation}

We  begin by proving tight upper bounds for the sum of $\frac{\log p}p$ over primes $p$ with $(d/p)=1$ using an old idea of Selberg, followed by 
proving tight upper bounds for the sum of $\frac{\log p}p$ over primes $p$ with $(d/p)=-1$ via an argument initiated by OpenAI.  Both of these arguments use only elementary ideas.   The upper bound over primes $p$ with $(d/p)=-1$ can be significantly improved by using the usual (though non-elementary) machinery taught in analytic number theory, and we do so as this will eventually lead to better upper bounds on $ (L'/L)(1,(d/\cdot))$ as well as better lower bounds for the distance from a real zero to the $1$-line.

      Based on the argument that Selberg \cite{selberg49,selberg} used to prove that $L(1,(d/\cdot))\ne 0$ we will establish the following:
 
\begin{theorem}\label{thm:Selberg}
Let $d$ be a fundamental discriminant.  For any integer $x\geq 3$   we have
\[
  \sum_{\substack{p\leq x \\ (d/p)=1}} \frac{ \log p}p   \leq \frac 12 \log  x+\frac 12 \log |d|+ O(1).
\]
\end{theorem}

We also develop a proof of a matching bound based on OpenAI's breakthrough \cite{openai}, as developed by Green \cite{green}, van Dijk \cite{vandijk}  and Pollack \cite{poll}   (who were influenced by preprints of Banks \cite{banks}, Barreto  \cite {barreto}, Kerr   \cite{kerr}; and see also the discussion in     Tao \cite{tao}):

\begin{theorem}\label{thm:NoSiegel}
Let $d$ be a fundamental discriminant.  For any integer $x\geq 3$   we have
 \[
   \sum_{\substack{p\leq x \\ (d/p)=-1}} \frac{\log p} p\leq \frac 12 \log x+\frac 14 \log |d|+\log\log x+O(1).
 \]
\end{theorem}

We can use these results to bound how close to $1$ a real zero $\beta$ of $L(s,(d/\cdot))$ can be, then use this bound to improve Theorem \ref{thm:NoSiegel} for large $x$, and so on; bootstrapping in this way leads to Corollaries \ref{thm:L'/L} and \ref{cor:L2}, while Corollary \ref{cor:L1} follows directly from Theorem \ref{thm:NoSiegel}.
We complete this section by describing the bootstrapping process, and then prove Theorems \ref{thm:Selberg} and \ref{thm:NoSiegel} in the following two sections.

Throughout this  section $\chi=(d/\cdot)$ and $|d|$ is sufficiently large. 

\subsection{Bounding the tail}
To control the tail of $\sum_p\frac{\chi(p)\log p}p$ we will prove  the following which uses standard (non-elementary) facts about Dirichlet $L$-functions.
 
 \begin{lemma}\label{lem:tail}
There is an absolute constant $C$ such that the following holds for $|d|$ sufficiently large.
Let $0<\kappa\le\frac12$ and suppose that $L(\beta,\chi)\neq 0$ for all real $\beta\in[1-\kappa,1)$. If
\[
\log x\ \geq\ C\log|d|\log\log|d| +\frac3\kappa\log\frac1\kappa
\]
then  $\Big|\sum_{x<p\le x'}\frac{\chi(p)\log p}p\Big|\ll \frac{1}{\log x}$ for all $x'>x$,
and so $\sum_p\frac{\chi(p)\log p}p$ converges.
\end{lemma}

\begin{proof}
Let $c_0\in(0,1]$ be the constant in the classical zero-free region \cite[Ch.~14]{dave}: $L(s,\chi)$ has at
most one zero $\beta_1$ (real and simple) with $\sigma\ge 1-c_0/\log(|d|(|t|+2))$, and by hypothesis
$\beta_1<1-\kappa$ if it exists. Since $\chi$ is real the non-trivial zeros are symmetric under
$\rho\mapsto1-\overline\rho$, so every non-trivial zero other than $\beta_1,1-\beta_1$ satisfies
$c_0/\log(|d|(|\gamma|+2))\le\beta\le1-c_0/\log(|d|(|\gamma|+2))$. For $|d|\le t$ and $2\le T\le t$ the explicit formula
\cite[Chs.~19--20]{dave} and $N(T+1,\chi)-N(T,\chi)\ll\log(|d|T)$ give
\[
\frac{|\psi(t,\chi)|}t\ \ll\ t^{-\kappa}+\frac{t^{-1/2}}\kappa+ \log^2(|d|T)
\exp\Big(-\frac{c_0\log t}{\log(|d|T)}\Big)+\frac{\log^2t}T ,
\]
the first two terms coming from $\beta_1$ and $1-\beta_1$.

If $C\log|d|\log\log |d| \leq \log t \leq (\log|d|)^2$ let $T=e^{\log t/\log|d|}$, so that $\log(|d|T)\le 2\log|d|$;
if $\log t > (\log|d|)^2$ let $T=e^{\sqrt{\log t}}$, so that $\log(|d|T)\le 2\sqrt{\log t}$. In either case, if $C$
is sufficiently large in terms of $c_0$, then
\[
\log^2(|d|T) \exp \Big(-\frac{c_0\log t}{\log(|d|T)} \Big)+\frac{\log^2t}T\ll \frac{1}{(\log t)^2}.
\]
Therefore, by partial summation, for $x<x'$,
\[
\sum_{x<n\le x'}\frac{\chi(n)\Lambda(n)}n=\int_x^{x'}\frac{d\psi(t,\chi)}t 
\ll  \frac{x^{-\kappa}}\kappa   + \frac{1}{\log x},
\]
and the prime powers contribute $O(x^{-1/2}\log x)$. Now $\kappa\log x\geq 3\log\frac1\kappa$, so  $\frac{x^{-\kappa}}\kappa\ll\frac1{\log x}$, and the result follows.
\end{proof}

\subsection{The bootstrapping process}

For $x\geq 2$, since $\sum_{p\le x}\frac{\log p}p=\log x+O(1)$, we have
\begin{equation}\label{eq:pm}
\begin{aligned}
\sum_{p\leq x}\frac{\chi(p)\log p}p&=2\sum_{\substack{p\leq x\\ \chi(p)=1}}\frac{\log p}p-\log x+E_d(x)+O(1)\\
&=\log x-2\sum_{\substack{p\leq x\\ \chi(p)=-1}}\frac{\log p}p-E_d(x)+O(1),
\end{aligned}
\end{equation}
where $E_d(x):=\sum_{p\le x,\ p\mid d}\frac{\log p}p\le E_d:=\sum_{p\mid d}\frac{\log p}p\le \log\log|d|+O(1)$. Moreover, since the prime powers contribute $O(1)$,
\begin{equation}\label{eq:L'/L}
\frac{L'}{L}(1,\chi)=-\lim_{x\to\infty}\sum_{p\leq x}\frac{\chi(p)\log p}p+O(1),
\end{equation}
the limit existing by Lemma \ref{lem:tail}.

\medskip\noindent\textbf{Step 1.} By Dirichlet's class number formula we have $L(1,\chi)\gg|d|^{-1/2}$; indeed for $d<0$ it equals $2\pi h(d)/(w_d\sqrt{|d|})\geq \pi/(3\sqrt{|d|})$, and for $d>0$ it equals $2h(d)\log\varepsilon_d/\sqrt d\geq 2\log(\frac{1+\sqrt5}2)\big/\sqrt d$. Partial summation and the P\'olya--Vinogradov inequality give $|L'(\sigma,\chi)|\ll\log^2|d|$ in the range $1-1/\log|d|\le\sigma\le 1$. Hence if $L(\beta,\chi)=0$ with $\beta$ real then, from the identity $L(1,\chi)=\int_\beta^1L'(\sigma,\chi)\,d\sigma$, we get $1-\beta\geq \kappa_0:=c_1/(\sqrt{|d|}\log^2|d|)$ for some absolute $c_1>0$. Lemma \ref{lem:tail} applies with $\kappa=\kappa_0$ and $x_0:=\exp(\sqrt{|d|}(\log|d|)^4)$, so for $x>x_0$, as no prime $p>x_0$ divides $d$,
\[
   \sum_{\substack{x_0<p\leq x \\ \chi(p)=-1}} \frac{\log p} p =\frac12\sum_{x_0<p\le x}\frac{\log p}p-\frac12\sum_{x_0<p\le x}\frac{\chi(p)\log p}p= \frac 12 \log (x/x_0)+O(1).
\]
Combining this with Theorem \ref{thm:NoSiegel} at $x_0$ gives, for all $x\geq 3$,
 \[
   \sum_{\substack{p\leq x \\ \chi(p)=-1}} \frac{\log p} p\leq \frac 12 \log x+\frac 14 \log |d|+\min\{ \log\log x,  \log\log x_0\}+O(1).
 \]
 By \eqref{eq:pm} and \eqref{eq:L'/L}, and since $\log\log x_0=\frac12\log|d|+4\log\log|d|$, this implies that
 \[
 \frac{L'}{L}(1,\chi)  = \lim_{x\to \infty}\bigg( 2\sum_{\substack{p\leq x \\ \chi(p)=-1}} \frac{\log p} p-\log x\bigg)+  E_d +O(1)
 \leq \frac 32  \log |d| + 9\log\log |d|+O(1).
  \]
 Now \cite[Ch.~12]{dave}  implies  that
\[
\frac{L'}{L}(1,\chi) \geq  \sum_{\substack{\rho: L( \rho, \chi)=0 \\ |\text{Im}(\rho)|<1}} \frac 1{1-\rho} -\frac12\log |d|+O(1).
\]
Now $\frac 1{1-\rho} + \frac 1{1-\overline{\rho} }$ is a positive real number for every non-real $\rho$ (and $\frac1{1-\rho}>0$ for real $\rho$), so if $L(\beta,\chi)=0$ with $\beta\in \mathbb R$ then
\begin{equation}\label{eq:betaL'}
\frac{L'}{L}(1,\chi) \geq \frac 1{1-\beta} -\frac12\log |d|+O(1).
\end{equation}
Combining these estimates gives
\[
\frac 1{1-\beta}   \leq  2  \log |d| + 9\log\log |d|+O(1)\ \le\ 3\log|d| .
\]

\medskip\noindent\textbf{Step 2.} We now apply Lemma \ref{lem:tail} with $\kappa=\kappa_1:=1/(3\log|d|)$. There is an absolute constant $C_1$ such that, with $x_1:=\exp(C_1\log|d|\log\log|d|)$, the argument of Step 1 gives, for all $x\ge3$,
 \begin{equation} \label{eq: L final}
   \sum_{\substack{p\leq x \\ \chi(p)=-1}} \frac{\log p} p\leq \frac 12 \log x+\frac 14 \log |d|+\min\{ \log\log x,  \log\log x_1\}+O(1),
  \end{equation}
where now $\log\log x_1=\log\log|d|+\log\log\log|d|+O(1)$.

\begin{proof} [Proofs of Corollaries \ref{thm:L'/L} and \ref{cor:L2}]
By \eqref{eq:pm}, \eqref{eq:L'/L} and \eqref{eq: L final},
\[
\frac{L'}{L}(1,\chi)\le \frac12\log|d|+2\log\log x_1+E_d+O(1)\le \frac12\log |d|+3\log\log |d|+2\log\log\log |d|+O(1),
\]
and so, by \eqref{eq:betaL'}, if $L(\beta,\chi)=0$ with $\beta$ real then
\[
\frac 1{1-\beta}   \leq   \log |d|+3\log\log |d| + 2\log\log\log |d|+O(1),
\]
which is Corollary \ref{cor:L2}. In the other direction, by \eqref{eq:pm} and Theorem \ref{thm:Selberg}, for every $x\ge 3$,
\[
   \sum_{ p\leq x } \frac{\chi(p)\log p} p= 2  \sum_{\substack{p\leq x\\ \chi(p)=1}} \frac{ \log p}p  -\log x+E_d(x)+O(1)
   \leq   \log |d|+E_d+ O(1) \leq \log |d| +\log\log |d|+O(1),
 \]
so that $-\frac{L'}{L}(1,\chi)\le\log|d|+\log\log|d|+O(1)$ by \eqref{eq:L'/L}. Corollary \ref{thm:L'/L} follows. All of the implicit constants are effective.
\end{proof}

\begin{proof} [Proof of Corollary \ref{cor:L1}]
Since $(1*\chi)(n)=\sum_{m\mid n}\chi(m)\ge 0$ and
$(1*\chi)(p)=2$ when $\chi(p)=1$, the hyperbola method and the P\'olya--Vinogradov
inequality $|\sum_{n\le N}\chi(n)|\le\sqrt{|d|}\log|d|$ give, for all $t\ge1$,
\[
\#\{p\le t:\ \chi(p)=1\} \le\ \frac12\sum_{n\le t}(1*\chi)(n)\ \le\ \frac12\,tL(1,\chi)+3\sqrt{t|d|}\log|d| .
\]
 By partial summation (using that $\frac{\log t}t$ is decreasing), for $x\geq y:=|d|(\log|d|)^5$ we get
\[
\sum_{\substack{y<p\le x\\ \chi(p)=1}}\frac{\log p}{p}
\le    \frac {L(1,\chi)}4\log^2x +O\Big(\frac{\sqrt{|d|}\log|d|\,\log y}{\sqrt y}\Big)=\frac {L(1,\chi)}4\log^2x +o(1).
\]
On the other hand no prime $p>y$ divides $d$, so for $x=|d|^B$ with $B\geq 2$ fixed, Theorem \ref{thm:NoSiegel} (with $\log\log x=\log\log|d|+O(1)$) gives
\[
\sum_{\substack{y<p\le x\\ \chi(p)=1}}\frac{\log p}{p}\geq \sum_{y<p\le x }\frac{\log p}{p}- \sum_{\substack{ p\le x\\ \chi(p)=-1}}\frac{\log p}{p}
\ge  \frac 12 \log x-\frac 54 \log |d|-6\log\log|d| +O(1).
\]
Comparing these two inequalities gives $L(1,\chi)\geq \frac{2B-5}{B^2\log|d|}\big(1-O(\frac{\log\log|d|}{\log|d|})\big)$, and $(2B-5)/B^2$ is maximised at $B=5$, where it equals $\frac15$. Note that this proof is entirely elementary.
\end{proof}

\section{Lattice points}

\begin{proof} [Proof of Theorem \ref{thm:Selberg}] We will assume that $x\geq |d|$, else the result follows from \eqref{eq: Trivial} immediately.
Let $H=(x|d|)^{1/2}$, so that $x\geq H\geq |d|$, and let $U:=\lfloor H\rfloor$ and $V:=\lfloor H/\sqrt{|d|}\rfloor=\lfloor\sqrt x\rfloor$.
We will consider the divisors of
\[
P=\prod_{\substack{1\leq u\le U\\ 1\le v\leq V}}  (u^2-dv^2).
\]
Now $P$ has $X:=UV=H^2/\sqrt{|d|}+O(H)$ factors, each non-zero (as $d\ne 1$ is not a square) and of absolute value $\leq 2H^2$, and so
\[
\log |P|  \leq X \log (2H^2) = X(\log x +\log |d|+O(1)).
\]

For an odd prime $p$ with $(d/p)=1$ fix $t$ with $t^2\equiv d \pmod p$ and set
\[
e_p:=\#\{ 1\leq u\le U,\ 1\le |v|\le V:\ u\equiv tv \pmod p\} .
\]
Each factor $u^2-dv^2$ of $P$ corresponds to the two pairs $(u,v)$ and $(u,-v)$. If $p\nmid v$ then at most one of $u\equiv tv$ and $u\equiv -tv \pmod p$ holds (as $p\nmid 2tv$), and if one does then $p\mid u^2-dv^2$. If $p\mid v$ and one of them holds then $p\mid u$, so both hold, and $p^2\mid u^2-dv^2$. Hence the number of the pairs $(u,\pm v)$ counted in $e_p$ is at most the exponent of $p$ in $u^2-dv^2$, and so $p^{e_p}$ divides $P$.
Therefore $\prod_{2<p\leq x,\ (d/p)=1} p^{e_p}$ divides $P$ and so
\[
\sum_{\substack{2<p\leq x\\ (d/p)=1}} e_p \log p\leq \log |P| \leq X (\log  x+\ \log |d|+O(1)).
\]

Now $\Lambda_{p,t}:=\{ (u,v)\in \Z^2: u\equiv tv \pmod p\}$ is a lattice of determinant $p$ containing $p\Z^2$. By periodicity, if $U_p,V_p$ are the least non-negative residues of $U,V \pmod p$ then
\begin{align*}
  \#\{ 1\leq u\leq U,\ 1\leq |v| \leq V:\ &u\equiv tv \pmod p\} - \frac {2UV}p \\
  &= \#\{ 1\leq u\leq U_p,\ 1\leq |v| \leq V_p:\ u\equiv tv \pmod p\} - \frac {2U_pV_p}p ,
\end{align*}
and this is $0$ if $U_p=0$ or $V_p=0$.

 Let $\lambda_p:=\max\{\min(U_p/|a|, V_p/|b|) :\ (a,b)\in\Lambda_{p,t}\setminus\{0\}\}$ (where $c/0:=\infty$). Rescaling by
$(u,v)\mapsto(u/U_p,v/V_p)$ maps the box to $(0,1]\times([-1,0)\cup(0,1])$ and $\Lambda_{p,t}$ to a lattice
of determinant $p/U_pV_p$ whose shortest non-zero vector (in the sup norm) has length $1/\lambda_p$. The standard estimate for the number of points of a lattice in a convex planar region (the area divided by the determinant, with an error that is $O(1)$ plus the perimeter divided by the length of the shortest non-zero lattice vector; see, e.g., \cite[Theorem~1.3]{bw}) therefore gives
\[
\Big|\#\{1\le u\le U_p,\ 1\le |v|\le V_p:\ u\equiv tv \ (p)\}-\frac{2U_pV_p}{p}\Big|\ \ll\  \lambda_p+1 .
\]
Now let $(a,b)\in\Lambda_{p,t}$ be any non-zero vector. Then $p\mid a^2-db^2\ne 0$ (as $d$ is not
a square), so $\max\{a^2,|d|b^2\}\ge p/2$; that is, $|a|\geq \sqrt{p/2}$ or $|b|\ge\sqrt{p/(2|d|)}$. Since $U_p,V_p<p$,
\[
\min\Big(\frac{U_p}{|a|},\frac{V_p}{|b|}\Big)\ \le\ \frac{p}{\max\{|a|,|b|\}}\ \le\  \sqrt{2|d|p}
\]
and since $U_p\le H$ and $V_p\le H/\sqrt{|d|}$,
\[
\min\Big(\frac{U_p}{|a|},\frac{V_p}{|b|}\Big)\ \le\  \frac{H}{\max\{|a|,\sqrt{|d|}\,|b|\}}\ \le\ \frac{H}{\sqrt{p/2}}.
\]
As these hold for every non-zero lattice vector, $\lambda_p\le\min\{\sqrt{2|d|p},\ \sqrt2\,H/\sqrt p\}$,
the first being the smaller when $p<H/\sqrt{|d|}$.
 Therefore
\[
\sum_{\substack{2<p\leq x\\ (d/p)=1}} e_p \log p\ \geq 2X \sum_{\substack{2<p\leq x\\ (d/p)=1}} \frac{ \log p}p - O \bigg( \sum_{p<H/\sqrt{|d|}}  \sqrt{|d|p}\log p+H\!\!\sum_{H/\sqrt{|d|}<p<x}  \frac{\log p} {\sqrt{p}}+\sum_{p\leq x} \log p\bigg),
\]
and this error term is $\ll   |d|^{-1/4} H^{3/2}+ Hx^{1/2}+x\ll X$, since $|d|^{-1/4} H^{3/2}=x^{3/4}|d|^{1/2}$ and $Hx^{1/2}=x|d|^{1/2}=X+O(H)$. Comparing the upper and lower bounds, dividing by $2X$, and noting that the prime $2$ contributes $O(1)$, the result follows.
   \end{proof}

\section{Proof of Theorem \ref{thm:NoSiegel}}

We may assume that $|d|$ is sufficiently large. If $x< 10 |d|^{1/2}$ the result follows from \eqref{eq: Trivial} immediately, so we assume that $x\ge 10|d|^{1/2}$. Let $q$ be the largest prime $\leq x$ with $(d/q)=-1$. If there is no such prime, or if $q<10|d|^{1/2}$, then the sum in Theorem \ref{thm:NoSiegel} is at most $\sum_{p\le 10|d|^{1/2}}\frac{\log p}p\le \frac12\log x+\frac14\log|d|+O(1)$. Otherwise $q\geq 10|d|^{1/2}$ and
\[
  \sum_{\substack{p\leq x \\ (d/p)=-1}} \frac{\log p} p=  \sum_{\substack{p\leq q \\ (d/p)=-1}} \frac{\log p} p ,
\]
so it suffices to prove the theorem with $x$ replaced by $q$. Henceforth we assume that $x=q$ is a prime with $(d/q)=-1$ and $q\geq 10|d|^{1/2}$. To avoid a clash with the parameter $y$ below, we write $x,z$ for the variables of our auxiliary polynomials (so from now on $x$ denotes a polynomial variable, and no longer the bound in Theorem \ref{thm:NoSiegel}).

 We define the set of algebraic integers
 \[
  \mathcal A:=\{ \alpha=u+v\sqrt{d}:\ 0<|u|\le H \text{ and } 0<|v|\le H/\sqrt{|d|}\}
 \]
where $H:=  q^{1/2}|d|^{1/4}$,
so that if $\alpha\in \mathcal A$ then $|\alpha|, |\overline{\alpha}|\leq 2H$.  Since $q\geq 10 |d|^{1/2}$  we have $6 |d|^{1/2}<2H<q$, and  
\[
|\mathcal A|=4\lfloor H\rfloor \Big\lfloor \frac{H}{\sqrt{|d|}}\Big\rfloor\ > 4 \cdot \frac{3H}4 \cdot  \frac{2H}{3\sqrt{|d|}} =2q .
\]
 
Let $y=\lfloor \log q \rfloor$ and $N_0=\lfloor q/y\rfloor$.  We will work with the set
\[
I:=\{ (i,j)\in \Z_{\geq 0}^2: i+yj<q\}.
\]
We will need some statistics about this set:

\begin{lemma} \label{lem: Size of I} Let $R:=\# I$, $D:=\max_{(i,j)\in I}  (i+jq)$ and $S:= \sum_{(i,j)\in I}    i$. Then, as $q\to\infty$,
\[
R\sim \frac{N_0q}2\sim \frac{q^2}{2y},\qquad D\sim N_0q\sim \frac{q^2}y \qquad\text{and}\qquad S \sim \frac{q^3}{6y};
\]
moreover $D<2R$.
\end{lemma}

\begin{proof}  Write $q=N_0y+r$ with $0\leq r\leq y-1$. Then 
\[
R=\sum_{j=0}^{N_0}  (q-jy)=\frac{(N_0+1)}2 (2q-y N_0)=\frac{(N_0+1)}2 (q+r);
\]
 whereas
$D=\max_{(i,j)\in I}  (i+jq)  \leq N_0q+r<2R$.  The proof of the asymptotic for $S$ is similar.
\end{proof}

 In particular, with $N:=|\mathcal A|\cdot N_0$ we have, since   $|\mathcal A|\ge 2q $,
\[
N-2R\ \ge\ 2qN_0-(N_0+1)(q+r)\ =\ N_0y(N_0-1)-2r\ >\  y(N_0^2-N_0-2) \geq  0
\]
as $r<y$ and $N_0\geq 2$.

If $p$ is an odd prime with $(d/p)=-1$ then
\[
\alpha^p = (u+v\sqrt{d})^p\equiv u^p+v^p d^{\frac{p-1}2} \sqrt{d}\equiv u-v\sqrt{d}=:\overline{\alpha} \pmod p
\]
 for every $\alpha=u+v\sqrt d\in\Z[\sqrt d]$. (All primes below are odd: $q$ is, and $p>y\geq 2$.)
If $q$ is a prime with $(d/q)=-1$ then the elements of $ \mathcal A$ are distinct mod $q$ provided $2H<q$  (if $u+v\sqrt d\equiv u'+v'\sqrt d \pmod q$ then $q\mid u-u'$ and $q\mid v-v'$, whereas $|u-u'|,|v-v'|\le 2H<q$). Moreover $\alpha^q\equiv \overline{\alpha}\pmod q$ for all $\alpha\in\mathcal A$, and $\Z[\sqrt{d}]/q\Z[\sqrt{d}] \cong \mathbb F_{q^2}$ is a field as $(d/q)=-1$.

\begin{lemma} \label{lem: Matrix}
Let $\mathcal N=\{ (\alpha,n): \alpha\in \mathcal A \text{ and } 0\leq n\leq N_0-1\}$ so that $N:=|\mathcal N|=|\mathcal A|\cdot N_0\geq 2R$  (as shown above).
 Define the  matrix $B$ with rows indexed by elements of $I$ and columns indexed by  elements of $\mathcal N$ where
\[
B_{(i,j),(\alpha, n)}:= \bigg(\frac{\partial}{\partial x}\bigg)^n x^iz^j\bigg|_{x=\alpha, z=\overline{\alpha}}
= \begin{cases}   \frac{i!}{(i-n)!}  \alpha^{i-n} \overline{\alpha}^j & \text{ if } n\leq i;\\   0 & \text{ otherwise.}    \end{cases}
\]
The rows of $B$ are linearly independent  (over $\mathbb C$).
\end{lemma}

\begin{proof} Suppose that  the rows are linearly dependent and so there are coefficients $c_{i,j}$, not all zero, for which
$ \sum_{(i,j)\in I} c_{i,j}\text{Row}_{i,j}(B)=0$.  The entries of $B$ lie in $K=\Q(\sqrt d)$, and the rank of a matrix does not change under field extension, so we may take the $c_{i,j}\in K$; multiplying through by a suitable integer we may assume $c_{i,j}\in \Z[\sqrt d]$, and then dividing by the largest power of $q$ dividing all of the $c_{i,j}$ (which is possible since $q\Z[\sqrt d]=\{u+v\sqrt d:\ q\mid u,\ q\mid v\}$) we may assume that not all $c_{i,j}\equiv 0 \pmod q$.  

Now, for all $(\alpha,n)\in \mathcal N$ we have
\[
\sum_{(i,j)\in I}   c_{i,j} B_{(i,j),(\alpha, n)}= \bigg(\frac{\partial}{\partial x}\bigg)^n C(x,z) \bigg|_{x=\alpha, z=\overline{\alpha}} =0
\text{ where } C(x,z):=\sum_{i,j} c_{i,j}  x^iz^j.
\]

Define $f(x)=C(x,x^q)=\sum_{(i,j)\in I} c_{i,j}  x^{i+jq}$.  Now if $(i,j)\in I$ then $0\leq i\leq q-1$ and  so the exponents $i+jq, (i,j)\in I$ are all distinct, which implies that  $f(x)\not\equiv 0\pmod q$ (as some $c_{i,j}\not\equiv 0\pmod q$) and that $\deg f\le D$. Since $\frac{\partial}{\partial x} x^q\equiv 0 \pmod q$ we deduce that
\[
f^{(n)}(x) \equiv \bigg(\frac{\partial}{\partial x}\bigg)^n C(x,z) \bigg|_{z=x^q} \pmod q
\]
for $0\leq n\leq N_0-1$; explicitly, comparing coefficients, $\frac{(i+jq)!}{(i+jq-n)!}=\prod_{m=0}^{n-1}(i+jq-m)\equiv \frac{i!}{(i-n)!}\pmod q$ if $n\le i$, while if $n>i$ the product contains the factor $jq$ and so is $\equiv 0\pmod q$.
Taking $x=\alpha$ so that $z=\alpha^q\equiv \overline{\alpha}\pmod q$, we deduce that  $ f^{(n)}(\alpha)\equiv 0 \pmod q$ for $0\le n\le N_0-1$.
Therefore $f(x) \pmod q$ has a zero of order $\geq N_0$ at each $\alpha\in\mathcal A$ (here we use $N_0\le q$: writing $f(x)=\sum_k b_k (x-\alpha)^k$ we have $f^{(n)}(\alpha)=n!\,b_n$, and $n!$ is invertible mod $q$ for $n<N_0\le q$, so $b_0\equiv\dots\equiv b_{N_0-1}\equiv 0\pmod q$). Since the $\alpha\in\mathcal A$ are distinct mod $q$ and $f\not\equiv 0\pmod q$ is a polynomial over the field $\mathbb F_{q^2}$,
\[
2R\leq  N=  |\mathcal A|\cdot N_0 \leq  \deg (f \bmod q)\le D  = \max_{(i,j)\in I} (i+jq)<2R
\]
by Lemma \ref{lem: Size of I}, a contradiction.
\end{proof}

 Lemma \ref{lem: Matrix}  implies that there is an $R$-by-$R$ submatrix $M$ of $B$ with non-zero determinant.
We order the rows of $M$ by increasing values of $i+jy$ (and if two are equal by increasing $i$-value).

\begin{lemma} \label{lem: prime div} [Divisibility by small primes]  For any given prime $p$ with $y<p< q$ and $(d/p)=-1$ we have $p^{e_p}$ divides $\det M$ in $\Z[\sqrt d]$, where
$e_p:=\sum_{(i,j)\in I}\lfloor i/p\rfloor >S/p-R$.
\end{lemma}

\begin{proof}  If $i=\ell p+r, 0\leq r\leq p-1$ then we replace the row indexed by $x^i z^j $ by the row corresponding to
\[
x^r z^j (x^p-z)^\ell = \sum_{k=0}^\ell   \binom{\ell}k (-1)^k x^{i-kp}z^{j+k},
\]
 which is $x^iz^j$  (the $k=0$ term) plus an integer linear combination of earlier rows  ---  the indices $(i-kp,j+k)$ with $1\le k\le \ell$ do lie in $I$, since $i-kp\ge 0$ and 
 \[
 (i-kp)+(j+k)y = i+jy -k(p-y)< i+jy
 \]
   as $p>y$, and they precede $(i,j)$ in our ordering. Doing this for every row, the matrix of the change of rows is unitriangular with respect to our ordering, and so this new matrix has the same determinant as the old.  We claim that all of the entries in this  new row are divisible by $p^\ell$ so we may divide that out:
   Now $\frac{\partial}{\partial x}  ( \frac{x^p-z}p  )^\ell=\ell x^{p-1}  ( \frac{x^p-z}p  )^{\ell-1}$ and so
   \[
    \bigg(\frac{\partial}{\partial x}\bigg)^{{n}} x^rz^j \bigg( \frac{x^p-z}p \bigg)^\ell \text{ is a polynomial in } x, z \text{ and } \frac{x^p-z}p  \text{ with integer coefficients.}
   \]
 Evaluating at $x=\alpha,\ z=\bar\alpha$ we have $\frac{\alpha^p-\overline{\alpha}}{p}\in\Z[\sqrt d]$ because $(d/p)=-1$, and so every entry $\big(\frac{\partial}{\partial x}\big)^n x^rz^j(x^p-z)^\ell\big|_{x=\alpha,z=\bar\alpha}$ of the new row lies in $p^\ell\Z[\sqrt d]$, as claimed. Since the determinant is multilinear in the rows,
we deduce that $p^{e_p}$ divides $\det M$ where
\[
e_p= \sum_{(i,j)\in I}   \lfloor  i/p \rfloor >\sum_{(i,j)\in I}  \Big(  \frac ip -1\Big)  = \frac{S}{p} -R.\qedhere
\]
 \end{proof}

\begin{proof}   [Proof of Theorem \ref{thm:NoSiegel}]
Each $\alpha\in \mathcal A$ satisfies $|\alpha|, |\overline{\alpha}|\leq W:=2H=2q^{1/2}|d|^{1/4}$ and so
\[
|B_{(i,j),(\alpha, n)}| \leq  \bigg|\frac{i!}{(i-n)!}  \alpha^{i-n} \overline{\alpha}^j \bigg|\leq q^n W^{i+j-n}\leq   W^{i+j+n}\leq   W^{i+q/y+N_0}
\]
since $i<q$ and so $\frac{i!}{(i-n)!} <q^n$, {$q\le W^2$, $j<q/y$ and $n<N_0$}. Therefore {by expanding the determinant,}
\[
  |\det M| \leq R!  \cdot W^S \cdot W^{(q/y+N_0)R}
\]
and so since $W<q,\ R<q^2$ and $q/y+N_0\ll q/\log q$, we deduce that
\[
\log |\det M|  \leq  S\log W   +  O (qR )
\]
and the same bound holds for $\log |\det\overline M|$, where $\overline M$ is the matrix whose entries are the conjugates of those of $M$, since $|\alpha|,|\bar\alpha|\le W$ is symmetric.
Now $\prod_{y<p<q, (d/p)=-1} p^{e_p}$ divides $\det M$ {in $\Z[\sqrt d]$ by Lemma \ref{lem: prime div}}, and similarly $\prod_{y<p<q, (d/p)=-1} p^{e_p}$ divides $\det  \overline{M}$ by conjugating.  Multiplying these together gives  $\prod_{y<p<q, (d/p)=-1} p^{2e_p}$ divides $\det M\cdot \det \overline{M}{=N_{K/\Q}(\det M)}$ which is a non-zero {rational} integer {by Lemma \ref{lem: Matrix}}.  Therefore, using Lemma \ref{lem: prime div},
 \[
 \sum_{\substack{y<p<q \\ (d/p)=-1}}   \bigg(\frac{S}{p} -R\bigg) \log p  \leq    \sum_{\substack{y<p<q \\ (d/p)=-1}}   e_p \log p  \leq  {\tfrac12\log\big|N_{K/\Q}(\det M)\big|\le}  S\log W   +  O (qR ) {.}
 \]
 Now using that ${R}\sum_{p<q} \log p=O({qR})$  and  $qR\ll S$   we deduce that
 \[
 \sum_{\substack{y<p<q \\ (d/p)=-1}} \frac{\log p} p\leq  \log W +O(1)  = \frac 12 \log q+\frac 14 \log {|d|}+O(1).
 \]
 Now adding $   \sum_{p\leq y} \frac{\log p} p{=\log y+O(1)\le \log\log q+O(1)}$ to both sides, and noting that the prime $p=q$ contributes only $\frac{\log q}q=o(1)$, we obtain the claimed result.
 \end{proof}

 \section{Linnik's Theorem}

\subsection{Selberg's identities}
In his elementary proof of the prime number theorem for arithmetic progressions, Selberg \cite{selberg} noted that 
    \begin{equation}\label{eq:identity}
    \Lambda(n)\log n+\sum_{uv=n}\Lambda(u)\Lambda(v)=\sum_{gm= n}\mu(g)\log^2 m
    \end{equation}
and was able to deduce that, for fixed $Q$,
\[
\sum_{\substack{n\leq x \\  n\equiv a \pmod Q}}  \Lambda(n)\log n+\sum_{\substack{uv\leq x \\  uv\equiv a \pmod Q}}\Lambda(u)\Lambda(v)=\frac{2x\log x}{\phi(Q)}+O(x).
\]
He also observed that 
 \[
\Lambda(n)\log^2 n-2\sum_{tuv=n}\Lambda(t)\Lambda(u)\Lambda(v)   =  \sum_{gm= n}\mu(g)\log^2 m \log n-2 \sum_{gmr= n}\mu(g)\Lambda(r)\log^2 m
\]
Here we might hope to obtain 
\[
  2 \sum_{\substack{pqr \leq x \\ pqr\equiv a \pmod Q}}  \frac{\log p\log q \log r}{pqr}  -\sum_{\substack{p\leq x \\ p\equiv a \pmod Q}}  \frac{\log^3 p}p 
 =o\bigg( \frac{ (\log x)^3}{\phi(Q)}\bigg)
\]
though we can only see how to do so when $\frac{\log x}{\log Q}\to \infty$. Pollack \cite{poll} did prove something like this restricting attention to integers with all their prime factors $>x^\alpha$. Using integers with exactly  three prime factors to find primes is an idea from Friedlander and Iwaniec and  reminds one of the proof of Linnik due to Matom\"aki,   Merikoski, and   Ter\"av\"ainen \cite{mmt}.

 \subsection{Products of three primes}
 The following is essentially given in Pollack \cite{poll}:
\begin{proposition}  \label{prop:lower}
Fix $0<\lambda<1$ and $1<A<B$ with $B-A \geq 2\sqrt{2}$, and let $\alpha=A/B$.  If
\[
  \bigg(1+\frac 2{B-A} \log\frac{B-1}{A-1}  \bigg)^3< 2(1- \lambda )^2
\]
then for $Q$ sufficiently large and all $(a,Q)=1$
\[
   \sum_{\substack{Q^A<p,q,r \leq x=Q^B \\ pqr\equiv a \pmod Q}}  \frac{\log p\log q \log r}{pqr}> \lambda
  \frac{ ((1-\alpha)\log x)^3}{\phi(Q)}.
\]
\end{proposition}

\begin{proof} Let $\tau= \frac 1{B-A} \log\frac{B-1}{A-1} $, so that $(1+2\tau)^3< 2(1- \lambda )^2$, and let
\[
B_\chi := \sum_{x^\alpha<p\leq x}  \frac{\chi(p) \log p}p \quad\text{ and }\quad b_\chi=B_{\chi}/B_{\chi_0},
\]
so that $B_{\chi_0}=(1-\alpha)\log x+O(1)$. The sum in the proposition equals $\frac1{\phi(Q)}\sum_{\chi}\overline{\chi}(a)B_\chi^3$, so if the conclusion is false then
 \[
  \sum_{\chi \pmod Q} \overline{\chi}(a) b_\chi^3\leq \lambda(1+O(\tfrac 1{\log x}))
  \]
  and so
\[
 \sum_{\chi \ne \chi_0} |b_\chi|^3\geq\bigg| \sum_{\chi \ne \chi_0} \overline{\chi}(a) b_\chi^3\bigg|  \geq  1- \lambda-o(1).
\]
Now, by the Brun--Titchmarsh theorem\footnote{ Montgomery and Vaughan \cite[Theorem~2]{MV1}
showed that $\pi(y+h;Q,c)-\pi(y;Q,c)\ \le\ \frac{2h}{\phi(Q)\,\log(h/Q)}$ for $h>Q$.} and partial summation (writing $t=Q^s$, so that $\frac{\log t}{\log(t/Q)}=\frac s{s-1}$), for every $(b,Q)=1$,
\[
S_b:= \sum_{\substack{x^\alpha<p\leq x \\ p\equiv b \pmod Q}}  \frac{\log p}p \leq  \frac{(2+o(1))\log Q}{\phi(Q)}\int_A^B\frac{s\,ds}{s-1}= \frac{2( (1-\alpha)(1+\tau)+o(1))\log x}{\phi(Q)} ,
\]
and so
\begin{align*}
 \frac 1{\phi(Q)} \sum_{\chi \pmod Q} |B_\chi|^2 &=\sum_{b \pmod Q} S_b^2 \le \max_b S_b\cdot\sum_{b \pmod Q} S_b\\
 &\le \frac{2( (1-\alpha)(1+\tau)+o(1))\log x}{\phi(Q)}\cdot B_{\chi_0} = \frac{2((1-\alpha)^2(1+\tau)+o(1)) (\log x)^2}{\phi(Q)} ,
\end{align*}
so that
\[
 \sum_{\chi \ne \chi_0} |b_\chi|^2= \sum_{\chi  } |b_\chi|^2-1\leq 1+2\tau+o(1).
\]
Since $(\frac12+\tau)(1+2\tau)^2=\frac12(1+2\tau)^3<(1-\lambda)^2$, the fixed quantity $\eta:=\frac{1-\lambda}{1+2\tau}-\sqrt{\frac12+\tau}$ is positive, and
\[
\sqrt{\tfrac 1{2}+\tau}+\eta-o(1)\ \le\ \frac{\sum_{\chi \ne \chi_0} |b_\chi|^3}{\sum_{\chi \ne \chi_0} |b_\chi|^2}\ \leq\ \max_{\chi \ne \chi_0} |b_\chi|.
\]
Select $\chi$ to be a character $ \ne \chi_0$ with $ |b_\chi|$ maximal. The character $\chi$ is real, else $\overline{\chi}\ne\chi$ and we obtain
\[
1+2\tau +2\eta-o(1)\le |b_\chi|^2 +  |b_{\overline{\chi}}|^2\leq 1+2\tau+o(1),
\]
 a contradiction. So $\chi$ is induced by $(d/\cdot)$ for some fundamental discriminant $d\ne 1$ with $|d|\le Q$. No prime in $(x^\alpha,x]=(Q^A,Q^B]$ divides $Q$, so $B_\chi=S_d(Q^B)-S_d(Q^A)$ where $S_d(X):=\sum_{p\le X}\frac{(d/p)\log p}{p}$. By \eqref{eq:pm} together with Theorems \ref{thm:Selberg} and \ref{thm:NoSiegel}, for $X=Q^A$ or $Q^B$,
 \[
 -\tfrac12\log|d|-O(\log\log Q)\ \le\ S_d(X)\ \le\ \log|d|+O(\log\log Q),
 \]
 and therefore $|B_\chi|\le \frac32\log Q+O(\log\log Q)$. Since $(1-\alpha)\log x=(B-A)\log Q$ we deduce that
\[
2\le (B-A)\sqrt{\tfrac 1{2}+\tau}\ <\ (B-A)|b_\chi|\ =\ (1+o(1))\frac{|B_\chi|}{\log Q}\ \le\ \frac32+o(1),
 \]
 also a contradiction. (The argument shows that the hypothesis $B-A\geq 2\sqrt2$ could be weakened to $B-A\geq 3/\sqrt2$.)
 \end{proof}

 \subsection{Primes and products of three primes  in a residue class}
 
 Throughout, $A=22$, $B=46$ and $\alpha=A/B=11/23$. Let $x=Q^{B}$, $(a,Q)=1$, and
all $o(\cdot)$ and $O(\cdot)$ refer to $Q\to\infty$.
 
 Our main goal in the rest of this section will be to prove the following:
 
 \begin{proposition}\label{prop:technical}
Let $\kappa\ge0$ and suppose that 
\[
P_a:=\sum_{\substack{x^{3\alpha}<p\le x^{3}\\ p\equiv a\ (\mathrm{mod}\ Q)}}\frac{\log p}{p}=\sum_{\substack{Q^{66}<p\le Q^{138}\\ p\equiv a\ (\mathrm{mod}\ Q)}}\frac{\log p}{p} \le\frac{\kappa \log x}{\phi(Q)}.
\] Then
\begin{equation}\label{eq:thm}
T_a:=\sum_{\substack{x^{\alpha}<p,q,r\le x\\ pqr\equiv a\ (\mathrm{mod}\ Q)}}\frac{\log p\,\log q\,\log r}{pqr}\ \le\ \bigl(0.2183\,(1-\alpha)^3+\kappa+o(1)\bigr)\frac{(\log x)^3}{\phi(Q)},
\end{equation}
where the $o(1)$ does not depend on $\kappa$.
\end{proposition}

\begin{proof}[Proof of Theorem~\ref{thm:primeslinnik}, assuming Proposition~\ref{prop:technical}]
We have $B-A=24\ge2\sqrt2$ and
\[
\Bigl(1+\frac{2}{24}\log\frac{45}{21}\Bigr)^3=(1.063512)^3=1.20289,\qquad 1-\sqrt{1.20289/2}=0.22447,
\]
so Proposition~\ref{prop:lower} applies with $\lambda=0.2244$ and gives $T_a>0.2244\,(1-\alpha)^3(\log x)^3/\phi(Q)$ for $Q$ large. Let $\kappa=P_a\,\phi(Q)/{\log x}$; then Proposition~\ref{prop:technical} gives $T_a\le(0.2183\,(1-\alpha)^3+\kappa+o(1))(\log x)^3/\phi(Q)$. Comparing the two bounds,
\[
\kappa\ \ge\ (0.2244-0.2183)(1-\alpha)^3-o(1)=0.0061\times0.142024-o(1)=0.00086-o(1),
\]
so $P_a\ge0.00085\,{\log x}/\phi(Q)=0.0391\log Q/\phi(Q)>\log (Q^{138})/(3600\phi(Q))$ for $Q$ large.
\end{proof}
 
 \medskip
The proof of Proposition~\ref{prop:technical}, which is given in Appendix~\ref{app:B} (written with Claude's help) uses two tools and nothing else: the Brun--Titchmarsh inequality, and the exact linear sieve with logarithmic weights of Matom\"aki, Merikoski and Ter\"av\"ainen \cite[Lemma 8.1]{mmt}. The exact sieve gives an inequality relating the primes $\equiv a$ in a range $(x^{\rho_1},x^{\rho_2}]$ to the \emph{balanced} products $p_1p_2p_3\equiv a$ in that range (Lemma~\ref{lem:sieve}): the products are bounded by the primes plus the sieve's main term, which is small ($\approx4\sigma/B$ per unit of $\sigma=\log(p_1p_2p_3)/\log x$) because the level of distribution of the class $a$ up to $y$ is $1-\log Q/\log y$, close to $1$. The passage from $\sum1/(p_1p_2p_3)$ to the logarithmically weighted $T_a$ costs a factor $\frac19\sigma^3$ by the AM--GM inequality. The only triples the sieve does not see are those near the two corners of the cube $(x^\alpha,x]^3$ -- where the cube is so thin that Brun--Titchmarsh does better -- and a thin set of unbalanced triples (``slivers''); both are handled by Brun--Titchmarsh.
   \appendix
   \section{Proof of Proposition \ref{prop:technical} (Claude)}\label{app:B}

\medskip\noindent\textbf{Notation.} For a prime $p$ put $u_p=\log p/{\log x}$. For a triple of primes we write $u=(u_p,u_q,u_r)$ and $\sigma=u_p+u_q+u_r=\log(pqr)/{\log x}$. The cube is $\cC=(\alpha,1]^3$: the triples in \eqref{eq:thm} are those with $u\in\cC$, and then $\sigma\in(3\alpha,3]$. We write $\ell=1-\alpha=12/23$ and $L=\log x$.

\begin{lemma}\label{lem:BTregion}
Let $\cR\subset\cC$ be a finite union of convex polytopes. Then
\[
\sum_{\substack{p,q,r\ \mathrm{distinct},\ u\in\cR\\ pqr\equiv a\ (\mathrm{mod}\ Q)}}\frac{\log p\log q\log r}{pqr}\ \le\ \bigl(A_\ast\Vol(\cR)+o(1)\bigr)\frac{(\log x)^3}{\phi(Q)}
\]
where $A_\ast=\frac{2A}{A-1}=\frac{44}{21}=2.0952\ldots $.
\end{lemma}

\begin{proof}
Writing $\cR$ as a finite union of convex polytopes with disjoint interiors (the common boundaries have measure zero and, by the argument below, contribute $o((\log x)^3/\phi(Q))$), we may assume that $\cR$ is a convex polytope. Fix primes $q,r$; the set of $u_p$ with $(u_p,u_q,u_r)\in\cR$ is then an interval with endpoints $z_1\le z_2$ (possibly empty), and the condition $pqr\equiv a\pmod Q$ becomes $p\equiv b:=a\overline{qr}\pmod Q$. Split $(x^{z_1},x^{z_2}]$ into intervals $(t,t(1+\eta)]$ with $\eta=1/\log x$, apply the Brun--Titchmarsh inequality (footnote to Proposition~\ref{prop:lower}) to each, and sum. Writing $t=Q^s$, so that $s\in(A,B]$ and $\frac{\log t}{\log (t/Q)}=\frac{s}{s-1}\le\frac A{A-1}$, this gives, uniformly in $b$, $z_1$ and $z_2$,
\[
\sum_{\substack{x^{z_1}<p\le x^{z_2}\\ p\equiv b\ (\mathrm{mod}\ Q)}}\frac{\log p}{p}\ \le\ \frac{2+o(1)}{\phi(Q)}\int_{x^{z_1}}^{x^{z_2}}\frac{\log t}{\log (t/Q)}\,\frac{dt}t+o\Big(\frac{\log x}{\phi(Q)}\Big)\ \le\ \bigl(A_\ast(z_2-z_1)+o(1)\bigr)\frac{{\log x}}{\phi(Q)}.
\]
Now multiply by $\frac{\log q\log r}{qr}$ and sum over $q$ and $r$. By Mertens' theorem, in the form
\[\sum_{x^{w_1}<q\le x^{w_2}}\frac{\log q}q=(w_2-w_1)\log x+O(1),\]
the result is $\frac{(\log x)^3}{\phi(Q)}$ times a Riemann sum for $\iint (z_2-z_1)\,du_q\,du_r=\Vol(\cR)$, up to $o(1)$. The lemma follows.
\end{proof}

\subsection{The sieve inequality}

Put
\begin{equation}\label{eq:m}
m(\sigma)=\frac{2}{\sigma\theta(\sigma)}\log\frac{\sigma B}{\sigma B-2},\qquad \theta(\sigma)=1-\frac{1}{\sigma B}\qquad(\sigma\ge3\alpha).
\end{equation}

\begin{lemma}\label{lem:sieve}
Let $3\alpha\le\rho_1<\rho_2\le3$ with $\rho_2<2\rho_1$. Then
\begin{equation}\label{eq:sieve}
\sum_{(p_1,p_2,p_3)}\frac{1}{p_1p_2p_3}\ \le\ \sum_{\substack{x^{\rho_1}<p\le x^{\rho_2}\\ p\equiv a\ (\mathrm{mod}\ Q)}}\frac1p+\frac{1}{\phi(Q)}\Bigl(\int_{\rho_1}^{\rho_2}m(\sigma)\,d\sigma+o(1)\Bigr),
\end{equation}
where the sum is over the triples of primes $p_1,p_2,p_3\in(x^\alpha,x]$ with $k=p_1p_2p_3$ satisfying
\[
k\equiv a\ (\mathrm{mod}\ Q),\quad x^{\rho_1}<k\le x^{\rho_2},\quad p_2<p_1<k^{1/2},\quad p_2<p_3\le Qp_2^2 .
\]
\end{lemma}

\begin{proof}
Fix a small $\varepsilon>0$ and put
\[
X=x^{\rho_2},\qquad \alpha_0=\frac{\rho_1}{\rho_2}\in\Bigl(\frac12,1\Bigr),\qquad \theta_\varepsilon(\rho)=1-\frac{1}{\rho\rho_2B}-\varepsilon\quad(\alpha_0\le\rho\le1).
\]
We apply \cite[Lemma~8.1]{mmt} to the sequence $a_n=\mathbf 1_{n\equiv a\ (Q)}$, $n\le X$, $(n,Q)=1$, with their $\varepsilon$ equal to ours, their $Y=Q$, and their density function $h(d)=\mathbf 1_{(d,Q)=1}$. Its hypothesis is that for each $y=X^\rho\in[X^{\alpha_0},X]$ the sequence $(a_n)_{n\le y}$ has level of distribution $\theta_\varepsilon(\rho)$ with dimension $1$ and size $y/Q$, in the sense of \cite[Definition~1.3]{mmt}; this is \cite[Lemma~11.1]{mmt}, whose proof is one line: for $d\le y^{\theta_\varepsilon(\rho)}=y^{1-\varepsilon}/Q$ coprime to $Q$, $\#\{n\le y/d:\ dn\equiv a\ (Q)\}=\frac{y}{dQ}+O(1)$, and the errors add up to $O(y^{1-\varepsilon}/Q)$, which is $o$ of $\frac yQ\prod_{p<y}(1-\frac1p)$. (The function $\theta_\varepsilon$ is increasing, and $\theta_\varepsilon(\alpha_0)=1-\frac1{\rho_1B}-\varepsilon>\frac12$, as the lemma requires.) The conclusion of \cite[Lemma~8.1]{mmt} immediately implies, by discarding non-negative terms, that
\begin{equation}\label{eq:mmt}
\sum_{\substack{X^{\alpha_0}<p\le X\\ p\equiv a}}\frac1p\ \ge\ -M+ \sum_{(p_1,p_2,p_3)}\frac{1}{p_1p_2p_3}\ +o\Bigl(\frac{\log X}{Q}\prod_{\substack{p<X^{\alpha_0}\\ p\nmid Q}}\Bigl(1-\frac1p\Bigr)\Bigr),
\end{equation}

where
\[
M=\frac{e^{\gamma}\log X}{Q}\int_{\alpha_0}^{1}\int_{\theta_\varepsilon(\rho)/2}^{1/2}\ \prod_{\substack{p<X^{\rho(\theta_\varepsilon(\rho)-\beta)/2}\\ p\nmid Q}}\Bigl(1-\frac1p\Bigr)\frac{d\beta}{\beta}\,d\rho,
\]
and where we have kept only part of the term $\widetilde S_2$ of \cite[Lemma~8.1]{mmt}, discarding the rest of $\widetilde S_2$ and all of the $\widetilde S_r$ with $r\geq 4$. Indeed $\widetilde S_2$ is the sum of $a_k/k$ over the $k\in(X^{\alpha_0},X]$ written as $k=p_1p_2s$, where all of the prime factors of $s$ are $>p_2$, with
\[
 p_2<p_1<k^{1/2}\quad\text{and}\quad p_1p_2^3\ \ge\ k^{(1+\varepsilon)\theta_\varepsilon(\log k/\log X)} .
\]
Our triples have $s=p_3>p_2$ prime. Moreover if $k=X^\rho$ then $k^{\theta_\varepsilon(\rho)}=k^{1-\varepsilon}/Q$, so that $k^{(1+\varepsilon)\theta_\varepsilon(\rho)}=k^{1-\varepsilon^2}/Q^{1+\varepsilon}\le k/Q$; therefore our condition $p_3\le Qp_2^2$, which is the same as $p_1p_2^3\ge k/Q$, implies that $p_1p_2^3\ge k^{(1+\varepsilon)\theta_\varepsilon(\rho)}$.
Hence 
\begin{equation}\label{eq:S2}
\sum_{(p_1,p_2,p_3)}\frac{1}{p_1p_2p_3}\ \le\ \sum_{\substack{x^{\rho_1}<p\le x^{\rho_2}\\ p\equiv a}}\frac1p+M+o\Bigl(\frac{\log X}{Q}\prod_{\substack{p<X^{\alpha_0}\\ p\nmid Q}}\Bigl(1-\frac1p\Bigr)\Bigr).
\end{equation}
It remains to evaluate $M$ and the error term. For every $z\ge2$,
\begin{equation}\label{eq:mertens}
\prod_{\substack{p<z\\ p\nmid Q}}\Bigl(1-\frac1p\Bigr)\le\prod_{p\mid Q}\Bigl(1-\frac1p\Bigr)^{-1}\prod_{p<z}\Bigl(1-\frac1p\Bigr)=\frac{Q}{\phi(Q)}\cdot\frac{e^{-\gamma}}{\log z}\Bigl(1+O\Bigl(\frac1{\log z}\Bigr)\Bigr)
\end{equation}
by Mertens' theorem. With $z=X^{\alpha_0}$ this shows that the error term in \eqref{eq:S2} is $o(1/\phi(Q))$. In the double integral, $\theta_\varepsilon(\rho)-\beta\ge\theta_\varepsilon(\alpha_0)-\frac12>0$, so $z=X^{\rho(\theta_\varepsilon(\rho)-\beta)/2}\ge X^{c}$ for a fixed $c>0$ and \eqref{eq:mertens} gives $\prod_{p<z,p\nmid Q}(1-\frac1p)\le\frac{Q}{\phi(Q)}\cdot\frac{2e^{-\gamma}}{\rho(\theta_\varepsilon(\rho)-\beta)\log X}(1+O(1/{\log x}))$. Therefore
\[
M\le\frac{1+O(1/{\log x})}{\phi(Q)}\int_{\alpha_0}^1\frac{2}{\rho}\int_{\theta_\varepsilon(\rho)/2}^{1/2}\frac{d\beta}{\beta(\theta_\varepsilon(\rho)-\beta)}\,d\rho=\frac{1+O(1/{\log x})}{\phi(Q)}\int_{\alpha_0}^1\frac{2}{\rho\,\theta_\varepsilon(\rho)}\log\frac{1}{2\theta_\varepsilon(\rho)-1}\,d\rho,
\]
because $\int_{\theta/2}^{1/2}\frac{d\beta}{\beta(\theta-\beta)}=\frac1\theta\Bigl[\log\frac{\beta}{\theta-\beta}\Bigr]_{\theta/2}^{1/2}=\frac1\theta\log\frac1{2\theta-1}$. Substitute $\sigma=\rho\rho_2$, so that $X^\rho=x^\sigma$, $d\rho/\rho=d\sigma/\sigma$ and $\theta_\varepsilon(\rho)=\theta(\sigma)-\varepsilon$ with $\theta(\sigma)$ as in \eqref{eq:m}: the integral becomes $\int_{\rho_1}^{\rho_2}\frac{2}{\sigma(\theta(\sigma)-\varepsilon)}\log\frac{1}{2\theta(\sigma)-1-2\varepsilon}\,d\sigma$, and since $2\theta(\sigma)-1=\frac{\sigma B-2}{\sigma B}$ this tends to $\int_{\rho_1}^{\rho_2}m(\sigma)\,d\sigma$ as $\varepsilon\to0$. As the left side of \eqref{eq:S2} does not depend on $\varepsilon$, the lemma follows.
\end{proof}

For the proof of Proposition~\ref{prop:technical} we need the following consequence. Let $\cK$ be the set of $k=pqr$ with $p<q<r$ primes in $(x^\alpha,x]$ and $k\equiv a\pmod Q$, and for $k\in\cK$ let $\nu(k)\in\{0,1,2\}$ be the number of orderings $(p_1,p_2,p_3)$ of $\{p,q,r\}$ with $p_2<p_1<k^{1/2}$ and $p_2<p_3\le Qp_2^2$. Then the sum in \eqref{eq:sieve} is $\sum_{k\in\cK,\,x^{\rho_1}<k\le x^{\rho_2}}\nu(k)/k$. Put $\sigma_k=\log k/L$.

\begin{lemma}\label{lem:riemann}
Suppose $P_a\le\kappa L/\phi(Q)$. Then for any $3\alpha\le\rho_1<\rho_2\le3$,
\[
\sum_{\substack{k\in\cK\\ x^{\rho_1}<k\le x^{\rho_2}}}\frac{\sigma_k^3\,\nu(k)}{k}\ \le\ \frac1{\phi(Q)}\Bigl(\int_{\rho_1}^{\rho_2}\sigma^3m(\sigma)\,d\sigma+9\kappa+o(1)\Bigr).
\]
\end{lemma}

\begin{proof}
Let $J\ge2$ and $\rho^{(j)}=\rho_1+j(\rho_2-\rho_1)/J$ for $0\le j\le J$. Then $\rho^{(j+1)}-\rho^{(j)}\le\frac{3-3\alpha}{2}<3\alpha\le\rho^{(j)}$, so Lemma~\ref{lem:sieve} applies to each range $(x^{\rho^{(j)}},x^{\rho^{(j+1)}}]$, and on that range $\sigma_k^3\le(\rho^{(j+1)})^3\le\sigma^3+27(\rho_2-\rho_1)/J$ for all $\sigma\in[\rho^{(j)},\rho^{(j+1)}]$. Hence
\[
\sum_{\substack{k\in\cK\\ x^{\rho_1}<k\le x^{\rho_2}}}\frac{\sigma_k^3\nu(k)}{k}\le\sum_{j=0}^{J-1}(\rho^{(j+1)})^3\Bigl[\frac1{\phi(Q)}\Bigl(\int_{\rho^{(j)}}^{\rho^{(j+1)}}m(\sigma)\,d\sigma+o(1)\Bigr)+\sum_{\substack{x^{\rho^{(j)}}<p\le x^{\rho^{(j+1)}}\\ p\equiv a}}\frac1p\Bigr].
\]
The first part is at most $\frac1{\phi(Q)}\bigl(\int_{\rho_1}^{\rho_2}\sigma^3m(\sigma)\,d\sigma+O(1/J)+o(1)\bigr)$. In the second, write $\frac1p=\frac{\log p}{p}\cdot\frac{1}{\sigma_pL}$ with $\sigma_p=\log p/L$; for $p$ in the $j$-th range $(\rho^{(j+1)})^3/\sigma_p\le\sigma_p^2(1+O(1/J))\le9(1+O(1/J))$, so the second part is at most
\[
\frac{9(1+O(1/J))}{L}\sum_{\substack{x^{3\alpha}<p\le x^3\\ p\equiv a}}\frac{\log p}{p}\le\frac{9\kappa}{\phi(Q)}\bigl(1+O(1/J)\bigr).
\]
Let $Q\to\infty$ and then $J\to\infty$.
\end{proof}

The description of $\nu$ is simple. Let $k=pqr\in\cK$ with $p<q<r$. Since $p_1,p_3>p_2$, necessarily $p_2=p$, and the two candidates are $(p_1,p_2,p_3)=(r,p,q)$, admissible iff $r<k^{1/2}$ and $q\le Qp^2$, and $(p_1,p_2,p_3)=(q,p,r)$, admissible iff $r\le Qp^2$ (the condition $q<k^{1/2}$ is automatic). In terms of $u$, with $s_1>s_2>s_3$ the coordinates $u_r>u_q>u_p$ and $\log Q/L=1/B$:
\begin{equation}\label{eq:nu}
\nu=\mathbf 1\bigl[s_1<s_2+s_3,\ s_2\le2s_3+\tfrac1B\bigr]+\mathbf 1\bigl[s_1\le2s_3+\tfrac1B\bigr].
\end{equation}
In particular $\nu=2$ unless $s_1\ge s_2+s_3$ or $s_1>2s_3+\frac1B$ (the condition $s_2>2s_3+\frac1B$ implies the latter).

\medskip\noindent\textbf{Proof of Proposition~\ref{prop:technical}.} Suppose $P_a\le\kappa L/\phi(Q)$. We bound $T_a$ by splitting it into four pieces.

\subsection{Splitting the cube}
Triples with two equal primes contribute at most
\[
3\sum_{p>x^\alpha}\frac{\log^2p}{p^2}\sum_{r\le x}\frac{\log r}{r}\ll L^2x^{-\alpha}
\]
to $T_a$, which is negligible; so we consider only triples of distinct primes, each $k\in\cK$ arising from six of them. Put
\[
\rho_1=3\alpha+0.12=1.554783\ldots,\qquad \rho_2=3-0.16=2.84,
\]
and
\[
\cR_1=\{u\in\cC:\sigma\le\rho_1\},\qquad \cR_2=\{u\in\cC:\sigma>\rho_2\},\qquad \cD=\{u\in\cC:\nu\le1\}.
\]
A triple of distinct primes with $u\in\cC$ lies in $\cR_1$, in $\cR_2$, in $\cD$, or has $\rho_1<\sigma\le\rho_2$ and $\nu=2$; accordingly
\begin{equation}\label{eq:split}
T_a\le\Sigma(\cR_1)+\Sigma(\cR_2)+\Sigma(\cD)+\Sigma_{\mathrm{mid}}+o\Bigl(\frac{(\log x)^3}{\phi(Q)}\Bigr),
\end{equation}
where $\Sigma(\cR)$ is the sum in Lemma~\ref{lem:BTregion} and $\Sigma_{\mathrm{mid}}$ the sum over triples of distinct primes with $\rho_1<\sigma\le\rho_2$ and $\nu=2$. We bound the four pieces in units of $(\log x)^3/\phi(Q)$.

\subsection{The middle range}
For a triple with product $k$, the AM--GM inequality gives
\[
\log p\log q\log r\le(\tfrac13\log k)^3=\sigma_k^3(\log x)^3/27 .
\]
Grouping the six triples with the same product and using $1=\nu(k)/2$ there,
\[
\Sigma_{\mathrm{mid}}=6\sum_{\substack{k\in\cK,\ \nu(k)=2\\ x^{\rho_1}<k\le x^{\rho_2}}}\frac{\log p\log q\log r}{k}\le\frac{6(\log x)^3}{27}\sum_{\substack{k\in\cK,\ \nu(k)=2\\ x^{\rho_1}<k\le x^{\rho_2}}}\frac{\sigma_k^3\,\nu(k)}{2k}
\le\frac{(\log x)^3}{9}\sum_{\substack{k\in\cK\\ x^{\rho_1}<k\le x^{\rho_2}}}\frac{\sigma_k^3\nu(k)}{k}.
\]
By Lemma~\ref{lem:riemann},
\begin{equation}\label{eq:mid}
\Sigma_{\mathrm{mid}}\le\Bigl(\frac19\int_{\rho_1}^{\rho_2}\sigma^3m(\sigma)\,d\sigma+\kappa+o(1)\Bigr)\frac{(\log x)^3}{\phi(Q)}.
\end{equation}
The integral has an elementary upper bound: since $\log\frac{y}{y-2}\le\frac{2}{y-2}$ for $y>2$,
\[
\frac{\sigma^3m(\sigma)}{9}\le\frac{\sigma^3}{9}\cdot\frac{4B}{(\sigma B-1)(\sigma B-2)}=\Phi'(\sigma),\quad
\Phi(\sigma)=\frac{4}{9B^3}\Bigl[\frac{y^2}{2}+3y-\log(y-1)+8\log(y-2)\Bigr]_{y=\sigma B},
\]
by the partial fractions $\frac{y^3}{(y-1)(y-2)}=y+3-\frac1{y-1}+\frac8{y-2}$. With $B=46$: $\frac{4}{9B^3}=4.5661\cdot10^{-6}$, $\rho_1B=71.52$, $\rho_2B=130.64$, and
\begin{align*}
\Bigl[\tfrac{y^2}2+3y-\log(y-1)+8\log(y-2)\Bigr]_{y=130.64}&=8533.40+391.92-4.865+38.857=8959.32,\\
\Bigl[\tfrac{y^2}2+3y-\log(y-1)+8\log(y-2)\Bigr]_{y=71.52}&=2557.56+214.56-4.256+33.933=2801.79,
\end{align*}
so
\begin{equation}\label{eq:midnum}
\frac19\int_{\rho_1}^{\rho_2}\sigma^3m(\sigma)\,d\sigma\le\Phi(\rho_2)-\Phi(\rho_1)=4.5661\cdot10^{-6}\times6157.53=0.028116 .
\end{equation}
(Numerically the integral is $0.027835$.)

\subsection{The corners}
After the translation $t_i=u_i-\alpha$, the region $\cR_1$ is the simplex $\{t_i>0,\ t_1+t_2+t_3\le0.12\}$, which lies in $(0,\ell]^3$ since $0.12<\ell$; its volume is $0.12^3/6$. Likewise $\cR_2$ is the simplex of side $0.16$ at the vertex $(1,1,1)$. By Lemma~\ref{lem:BTregion},
\begin{equation}\label{eq:corners}
\Sigma(\cR_1)+\Sigma(\cR_2)\le\Bigl(\frac{A_\ast}{6}\bigl(0.12^3+0.16^3\bigr)+o(1)\Bigr)\frac{(\log x)^3}{\phi(Q)}=\bigl(0.002034+o(1)\bigr)\frac{(\log x)^3}{\phi(Q)}.
\end{equation}
(The sieve could be used on the corner ranges too, but there it is worse than Brun--Titchmarsh: the budget $\frac19\sigma^3m(\sigma)$ per unit of $\sigma$ does not vanish at the ends of the range, while the cross-section of the cube does, and the choice $0.12$, $0.16$ is roughly where the two cross.)

\subsection{The slivers}
By \eqref{eq:nu}, $\cD\subset\cD_1\cup\cD_2$ with $\cD_1=\{u\in\cC:s_1\ge s_2+s_3\}$ and $\cD_2=\{u\in\cC:s_1>2s_3+\frac1B\}$, where $s_1\ge s_2\ge s_3$ are the coordinates of $u$ in decreasing order; both are finite unions of convex polytopes, with volumes six times those of the sets of sorted triples $1\ge s_1\ge s_2\ge s_3>\alpha$. For $\cD_1$, with $t_i=s_i-\alpha$ and $1-2\alpha=\frac1{23}$,
\[
\Vol(\cD_1)=6\int\!\!\int_{\substack{\alpha<s_3\le s_2\\ s_2+s_3\le1}}(1-s_2-s_3)\,ds_2\,ds_3=3\int\!\!\int_{\substack{t_2,t_3\ge0\\ t_2+t_3\le1/23}}\Bigl(\frac1{23}-t_2-t_3\Bigr)dt_2\,dt_3
\]
\[
=\frac{1}{2\cdot23^3}=0.0000411 .
\]
For $\cD_2$: given $s_3$, the constraints are $2s_3+\frac1B<s_1\le1$ and $s_3\le s_2\le s_1$, which force $s_3<\frac12(1-\frac1B)=\frac{45}{92}$, while $s_3>\alpha=\frac{44}{92}$; so
\[
\Vol(\cD_2)=6\int_{44/92}^{45/92}\int_{2s_3+1/46}^{1}(s_1-s_3)\,ds_1\,ds_3=6\int_{44/92}^{45/92}\frac{(1-s_3)^2-(s_3+\frac1{46})^2}{2}\,ds_3
\]
\[
=3\cdot\frac{47}{46}\int_{44/92}^{45/92}\Bigl(1-2s_3-\frac1{46}\Bigr)ds_3=\frac{3\cdot47}{46\cdot92^2}=0.0003622,
\]
using $(1-s_3)^2-(s_3+\frac1{46})^2=(1-2s_3-\frac1{46})(1+\frac1{46})$ and $\int_{44/92}^{45/92}(\frac{90}{92}-2s_3)\,ds_3=\frac{1}{92^2}$. Lemma~\ref{lem:BTregion} gives
\begin{equation}\label{eq:slivers}
\Sigma(\cD)\le\bigl(A_\ast(0.0000411+0.0003622)+o(1)\bigr)\frac{(\log x)^3}{\phi(Q)}=\bigl(0.000845+o(1)\bigr)\frac{(\log x)^3}{\phi(Q)}.
\end{equation}

\subsection{Conclusion}
By \eqref{eq:split}, \eqref{eq:mid}, \eqref{eq:midnum}, \eqref{eq:corners} and \eqref{eq:slivers},
\[
T_a\le\bigl(0.002034+0.028116+0.000845+\kappa+o(1)\bigr)\frac{(\log x)^3}{\phi(Q)}=\bigl(0.030995+\kappa+o(1)\bigr)\frac{(\log x)^3}{\phi(Q)},
\]
and $0.030995/\ell^3=0.030995/0.142024=0.21824<0.2183$. This proves Proposition~\ref{prop:technical}. \qed


\end{document}